\documentclass[a4paper,12pt]{article}

\usepackage{amsmath,amssymb,latexsym,amsthm,amsfonts,setspace}
\usepackage[dvips]{color,graphicx}
\usepackage{float}

\theoremstyle{definition}
\newtheorem{definition}{Definition}[section]
\newtheorem{thm}[definition]{Theorem}

\newtheorem{cor}[definition]{Corollary}
\newtheorem{lem}[definition]{Lemma}
\newtheorem{exmp}[definition]{Example}
\newtheorem{rem}[definition]{Remark}

\newtheorem{question}{Question}

\theoremstyle{definition}
\numberwithin{equation}{section}

\begin{document}
\title{Bounds for the boxicity of Mycielski graphs}
\author{\footnote{Interdisciplinary Faculty of Science and Engineering, Shimane University, Shimane 
690-8504, Japan. \endgraf 
\hspace{0.27cm}{\it E-mail address}:\,kamibeppu@riko.shimane-u.ac.jp \endgraf
This work was supported by Grant-in-Aid for Young Scientists (B), No.25800091.}\,\,Akira Kamibeppu}
\date{}
\maketitle
\setstretch{1.2}

\begin{abstract}
A {\it box} in Euclidean $k$-space is the Cartesian product $I_1\times I_2\times \cdots \times I_k$, 
where $I_j$ is a closed interval on the real line.\ The {\it boxicity} of a graph $G$, denoted by $\text{box}(G)$, 
is the minimum nonnegative integer $k$ such that $G$ can be isomorphic to the intersection graph of a family of 
boxes in Euclidean $k$-space.\\ 
\indent Mycielski \cite{My55} introduced an interesting graph operation that extends a graph $G$ to a new 
graph $M(G)$, called the {\it Mycielski graph} of $G$.\ 
In this paper, we observe behavior of the boxicity of Mycielski graphs.\ 
The inequality $\text{box}(M(G))\geq \text{box}(G)$ holds for a graph $G$, and hence 
we are interested in whether the boxicity of the Mycielski graph of $G$ is more than that of $G$ or not.\ 
Here we give bounds for the boxicity of Mycielski graphs:\ for a graph $G$ with $l$ universal vertices, 
the inequalities $\text{box}(G)+\left \lceil \frac{l}{2}\right \rceil \leq \text{box}(M(G))\leq 
\theta (\overline{G})+\left\lceil \frac{l}{2}\right\rceil +1$ hold, where $\theta (\overline{G})$ is the 
{\it edge clique cover number} of the complement $\overline{G}$.\ Further observations determine the boxicity 
of the Mycielski graph $M(G)$, if $G$ has no universal vertices or odd universal vertices and satisfies  
$\text{box}(G)=\theta (\overline{G})$.\\ 
\indent We also present relations between the Mycielski graph $M(G)$ and its analogous ones $M_3(G)$ and $M_r(G)$ 
in the context of boxicity, which will encourage us to calculate the boxicity of $M(G)$ or $M_3(G)$. \\

\end{abstract}

\section{Introduction}
The notion of {\it boxicity} of graphs was introduced by Roberts \cite{Ro69}.\ 
It has applications in some research fields, like niche overlap in ecology (see \cite{Ro76, Ro78}) and 
fleet maintenance in operations research (see \cite{OR81}).\ Roberts \cite{Ro69} proved that the maximum 
boxicity of graphs with $n$ vertices is $\lfloor \frac{n}{2}\rfloor $ (also see \cite{CR83}), where 
$\lfloor x\rfloor $ denotes the largest integer at most $x$.\ Cozzens \cite{C81} proved that the 
task of computing boxicity of graphs is NP-hard.\ Some researchers have attempted to calculate or 
bound boxicity of graphs with special structure.\ 
Roberts \cite{Ro69} showed that the boxicity of a complete $k$-partite graph $K_{n_1, n_2, \ldots, n_k}$ 
is the number of $n_i$ which is at least 2.\ Scheinerman \cite{Sch84} proved that the boxicity of outer 
planar graphs is at most 2.\ Thomassen \cite{Tho86} proved that the boxicity of planar graphs is at 
most 3.\ Cozzens and Roberts \cite{CR83} investigated the boxicity of split graphs.\ As Chandran et al.\ 
\cite{CS07} say, not much is known about boxicity of most of the well-known graph classes.\ 
They proved that the boxicity of a graph $G$ is at most $\text{tw}(G)+2$, where $\text{tw}(G)$ is the 
treewidth of $G$, and presented upper bounds for chordal graphs, circular 
arc graphs, AT-free graphs, co-comparability graphs, and permutation graphs.\ Recently, 
Chandran et al.\ \cite{CDS09} found the following relation between boxicity and chromatic number.\  
\begin{thm}[\cite{CDS09}, Theorem 6.1]
Let $G$ be a graph with $n$ vertices.\ If $\text{box}(G)=\frac{n}{2}-s$ for $s\geq 0$, the 
inequality $\chi (G)\geq \frac{n}{2s+2}$ holds, where $\chi (G)$ is the {\it chromatic number} of $G$.\ 
\end{thm}
Theorem 1.1 implies that, if the boxicity of a graph with $n$ vertices is very close to the maximum 
boxicity $\lfloor \frac{n}{2}\rfloor $, the chromatic number of the graph must be very large.\ 
The converse does not hold in general; there is a graph whose boxicity is small, even if the 
chromatic number of the graph is large, like a complete graph.\ Also there are bipartite graphs 
with arbitrary large boxicity (see section 5.1 in \cite{CDS09} and also see \cite{CFM11}).\ 
However, a graph operation increasing chromatic number may admit increasing boxicity.\  
For example, the {\it join} of two graphs, taking the disjoint union of two graphs and 
adding all edges between them is desired one.\ Behavior of boxicity has been studied in the context 
of various graph operations (see \cite{CMS11, Tro79} for example).\ 
This paper is another attempt in this direction that studies behavior of boxicity in the context 
of Mycielski's graph operation.\ 

One of the purpose of this paper is to consider whether behavior of boxicity is similar to 
that of chromatic number under Mycielski's graph operation.\ Mycielski \cite{My55} invented an 
interesting graph operation that extends a graph $G$ to a new graph $M(G)$, 
called the {\it Mycielski graph} of $G$ or the {\it Mycielskian} of $G$.\ 
It is well-known that the chromatic number of the Mycielski graph of $G$ is more than 
that of $G$, actually, $\chi (M(G))=\chi (G)+1$ holds.\ We can construct (triangle-free) graphs 
with arbitrary large chromatic number by using the graph operation.\ 
Here we present the definition of the graph $M(G)$.\ Let $V(G)_i$ be a copy of the vertex set 
$V(G)$ of a graph $G$, where $i=1,2$.\ For each vertex $v\in V(G)$, the symbol $v_i$ denotes the vertex in $V(G)_i$ 
corresponding to $v$.\ The vertex set of $M(G)$ is defined to be $\{z\} \cup V(G)_1 \cup V(G)_2$, 
the disjoint union of the set of a single new vertex $z$ and copies $V(G)_1$ and $V(G)_2$, 
and the edge set of $M(G)$ is defined to be the union $E_1\cup E_2\cup E_3$, where
\begin{equation*}
E_1=\{u_1v_1\,|\,uv\in E(G)\},\, E_2=\{u_1v_2, v_1u_2\, |\, uv\in E(G)\},\, \text{and}\,\, E_3=\{zu_2 \, |\,u \in V(G)\}
\end{equation*}
and $E(G)$ denotes the edge set of $G$ (see Fig.\ 1 for example).\ 
Note that the inequality $\text{box}(M(G))\geq \text{box}(G)$ holds for a graph $G$ since $M(G)$ 
contains the subgraph induced by $V(G)_1$, isomorphic to $G$.\ 
\begin{figure}[!h]
\centering
\includegraphics[scale=1,clip]{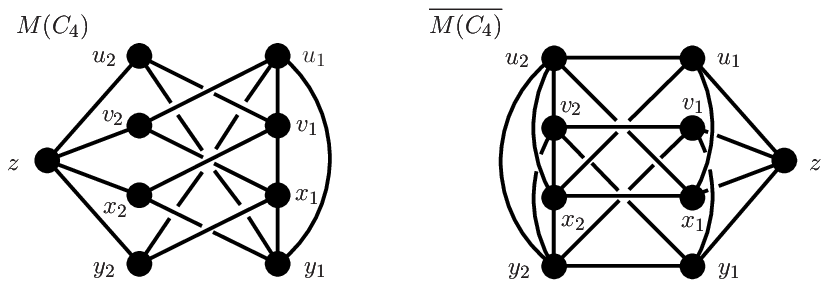}
\renewcommand{\baselinestretch}{1.3}
\caption[]{The Mycielski graph $M(C_4)$ of a cycle $C_4$ and its complement $\overline{M(C_4)}$.}
\end{figure}
So, first we are interested in whether 
the boxicity of the Mycielski graph $M(G)$ is more than that of $G$, 
the same as behavior of chromatic number under the graph operation, as mentioned at the beginning of this paragraph.\ 
Many researchers have studied Mycielski graphs and have compared a graph $G$ with $M(G)$ under 
various graph invariants (see \cite{LPU95, LWLG} for example).\   

In section 3, we improve the trivial lower bound for the boxicity of the Mycielskian of a graph $G$ in terms of 
the number of universal vertices of $G$.\ This implies that the boxicity of the Mycielski graph $M(G)$ is more 
than that of $G$ if the graph $G$ has universal vertices.\ Also note that there is 
a graph $G$ without universal vertices such that the boxicity of the Mycielski graph $M(G)$ is more than 
that of $G$.\ While such examples of graphs appear, there is also a graph $G$ such that 
$\text{box}(M(G))=\text{box}(G)$.\ 
As a conclusion, behavior of boxicity is not similar to that of chromatic number under Mycielski's graph operation 
in general.\ We reach the next purpose: Classify as many graphs as possible into $\text{box}(M(G))>\text{box}(G)$ 
or $\text{box}(M(G))=\text{box}(G)$.\

In section 4, we discuss upper bounds for the boxicity of Mycielski graphs.\ Chandran et al.\ 
\cite{CDS09} proved that the inequality $\text{box}(G)\leq \lfloor \frac{t(G)}{2}\rfloor +1$ holds for a graph $G$, 
where $t(G)$ is the minimum cardinality of a vertex cover of $G$.\ It is easy to see that 
$t(M(G))\leq 2t(G)+1$ for a graph $G$, and hence we have 
$\text{box}(M(G))\leq \lfloor \frac{t(M(G))}{2}\rfloor +1\leq t(G)+1$.\ Here we present another upper bound 
for the boxicity of the Mycielskian of a graph $G$ in terms of the edge clique cover number $\theta (\overline{G})$ 
of the complement $\overline{G}$.\ We also consider graphs that satisfy the equality 
$\text{box}(G)=\theta (\overline{G})$.\ The family of graphs satisfying $\text{box}(G)=\theta (\overline{G})$ 
contains complete multi-partite graphs, for example.\ Other examples of such graphs appear at the end of section 4.\ 
As a result, our observations determine the boxicity of their Mycielski graphs if original graphs have no universal 
vertices or odd universal vertices.\

In section 5, we consider relations between the Mycielski graph and its analogous one
$M_r(G)$, called the {\it generalized Mycielski graph} of $G$, in the context of boxicity, 
where $r\geq 3$.\ We present upper bounds for the boxicity of the generalized Mycielski graph $M_r(G)$ 
in terms of that of $M(G)$ for a bipartite graph $G$ or in terms of that of $M_3(G)$ for a graph $G$.\ 
These results will become our motivation a bit to calculate the boxicity of $M(G)$ or $M_3(G)$.\ 
 
\section{Preliminary}
\indent In this paper, all graphs are finite, simple and undirected.\ We use $V(G)$ for the vertex set 
of a graph $G$.\ We use $E(G)$ for the edge set of a graph $G$.\ An edge of a graph with endpoints $u$ 
and $v$ is denoted by $uv$.\ A vertex $v$ of $G$ is said to be {\it universal} if $v$ is adjacent to all 
vertices in $V(G)\setminus \{v\}$.\ A graph is said to be {\it trivial} if $E(G)$ is empty.\ For a subset 
$V$ of $V(G)$, let $G-V$ be the subgraph induced by $V(G)\setminus V$.\ For a subset $E$ of $E(G)$, let 
$G-E$ be the subgraph on $V(G)$ with $E(G)\setminus E$ as its edge set.\ 
A subset of $V(G)$ that induces a complete subgraph of $G$ is called a {\it clique} of $G$.\ 
For a graph $G$, its complement is denoted by $\overline{G}$.\ The {\it intersection graph} of a nonempty 
family $\mathcal{F}$ of sets is the graph whose vertex set is $\mathcal{F}$ and $F_1$ is adjacent to $F_2$ 
if and only if $F_1\cap F_2\ne \emptyset $ for $F_1$, $F_2\in \mathcal{F}$.\ The intersection 
graph of a family of closed intervals on the real line is called an {\it interval graph}.\ 
A graph $G$ can be represented as the intersection graph of a family $\mathcal{F}$ if there is a 
bijection between $V(G)$ and $\mathcal{F}$ such that two vertices of $G$ are adjacent if and only if 
the corresponding sets in $\mathcal{F}$ have nonempty intersection.\ 
A {\it box} in Euclidean $k$-space is the Cartesian product 
$I_1\times I_2\times \cdots \times I_k$, where $I_j$ is a closed interval on the real line.\ The 
{\it boxicity} of a graph $G$, denoted by $\text{box}(G)$, is the minimum nonnegative integer $k$ 
such that $G$ can be represented as (isomorphic to) the intersection graph of a family of boxes in 
Euclidean $k$-space.\ The boxicity of a complete graph is defined to be 0.\ If $G$ is an interval 
graph, $\text{box}(G)\leq 1$.\ If $H$ is an induced subgraph of $G$, $\text{box}(H)\leq \text{box}(G)$ 
holds by the definition.\ 

A graph is a {\it cointerval graph} if its complement is an interval graph.\ Lekkerkerker and 
Boland \cite{LB62} presented the forbidden subgraph characterization of interval or cointerval 
graphs.\ Cointerval graphs do not contain the complement of a cycle of length at least 4 as an induced 
subgraph, for example.\ It is easy to see that the union of a cointerval graph and isolated vertices 
is also a cointerval graph.\ A {\it cointerval edge covering} of a graph $G$ is a family $\mathcal{C}$ 
of cointerval spanning subgraphs of $G$ such that each edge of $G$ is in some graph of $\mathcal{C}$.\ 
For a set $X$, the cardinality of $X$ is denoted by $|X|$.\ The following theorem is useful to calculate of 
the boxicity of graphs.\
\begin{thm}[\cite{CR83}, Theorem 3]Let $G$ be a graph.\ Then, $\text{box}(G)\leq k$ if and only if 
there is a cointerval edge covering ${\mathcal C}$ of $\overline{G}$ with $|\,{\mathcal C}\,|=k$.
\end{thm}

\section{A lower bound for the boxicity of Mycielski graphs}
For a complete graph $K_n$, it is easy to see that $\text{box}(M(K_n))\geq 1>0=\text{box}(K_n)$ 
since $M(K_n)$ is not complete by the definition.\ We determine the boxicity of $M(K_n)$ 
next section (see Lemma 4.1).\ 
First we consider if the boxicity of the Mycielski graph of a graph $G$ is more than 
that of $G$ in general.\

\begin{question}
For a graph $G$, does the inequality $\text{box}(M(G))> \text{box}(G)$ hold?\
\end{question}
The following example shows that there exists a graph $G$ such that the equality 
$\text{box}(M(G))=\text{box}(G)$ holds.\ Here $C_n$ denotes a cycle with $n$ vertices.\ 
\begin{exmp}
The boxicity of the Mycielski graph of a cycle $C_4$ is equal to 2.\ To check this, we give a cointerval 
edge covering of the complement $\overline{M(C_4)}$ (see Fig.\ 1).\ 

Let $H_1$ and $H_2$ be the graphs appeared in Fig.\ 2.\ Both graphs are cointerval spanning subgraphs 
of $\overline{M(C_4)}$.\ 
\begin{figure}[!h]
\vspace{0.2cm}
\centering
\includegraphics[scale=1,clip]{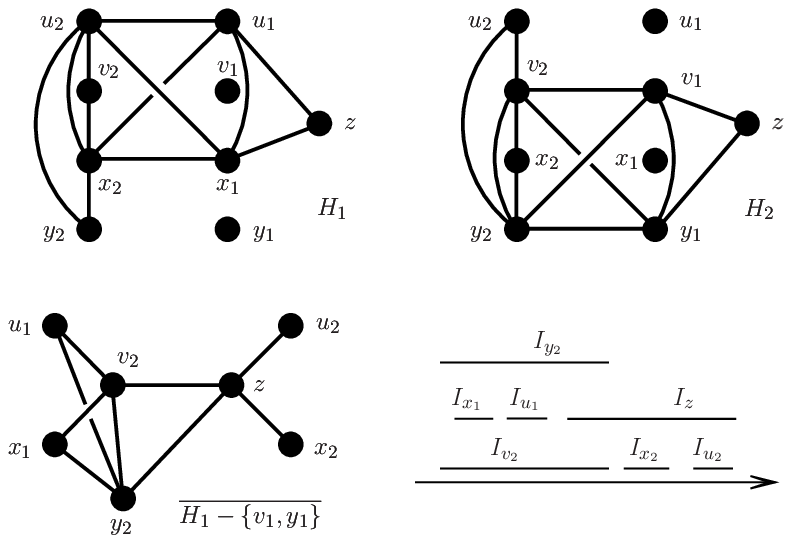}
\renewcommand{\baselinestretch}{1.3}
\caption[]{\parbox[t]{133mm}{Cointerval spanning subgraphs $H_1$ and $H_2$ of $\overline{M(C_4)}$ 
and an interval representation of $\overline{H_1-\{v_1, y_1\}}$.}}
\end{figure}
Note that the disjoint union of a cointerval graph and isolated 
vertices is also cointerval since these isolated vertices become pairwise adjacent universal vertices 
in the complement.\ 
Hence, we may prove that $H_1-\{v_1, y_1\}$ and $H_2-\{u_1, x_1\}$ are 
cointerval, instead of $H_1$ and $H_2$, respectively.\ A family of intervals on the real line 
with intersection graph isomorphic to $\overline{H_1-\{v_1, y_1\}}$ can be found as in the bottom 
of Fig.\ 2.\ 
Similar arguments work for $H_2-\{u_1, x_1\}$.\ Also see that $H_1$ 
and $H_2$ cover all edges of $\overline{M(C_4)}$.\ The family $\{H_1, H_2\}$ is 
a desired cointerval edge covering of $\overline{M(C_4)}$, and hence, $\text{box}(M(C_4))\leq 2$ 
by Theorem 2.1.\ Also note that $\text{box}(M(C_4))\geq \text{box}(C_4)=2$.
\end{exmp}
\begin{question}
Is there a graph $G$ such that the inequality $\text{box}(M(G))>\text{box}(G)$ holds?\ 
\end{question} 
\indent The {\it distance between two vertices $u$ and $v$} in a graph $G$ is defined by 
length of the shortest path from $u$ to $v$ in $G$ and is denoted by $d_G(u,v)$.\ If there 
exist no paths from $u$ to $v$ in $G$, define $d_G(u,v)=\infty $.\ Let $H_1$ and $H_2$ be subgraphs of $G$.\ 
The {\it distance between two subgraphs $H_1$ and $H_2$} in $G$, denoted by 
$d_G(H_1,H_2)$, is defined to be the minimum distance 
$\min \{d_G(v_1,v_2)\,|\, v_1\in V(H_1),\, v_2\in V(H_2)\}$.\ The following lemma is a generalization of Corollary 3.6 
in \cite{CR83}.\ 
\begin{lem}
Let $G$ be a graph and $H_1$, $H_2$ induced subgraphs of the complement $\overline{G}$.\ If 
$d_{\overline{G}}(H_1,H_2)\geq 2$, the following inequality holds: 
\begin{equation*}
\text{box}(G)\geq \text{box}(\overline{H_1})+\text{box}(\overline{H_2}).
\end{equation*} 
\begin{proof}
If either $H_1$ or $H_2$ is trivial, say $H_1$, then $\overline{H_1}$ is complete.\ Hence, 
$\text{box}(\overline{H_1})=0$.\ Since $\overline{H_2}$ is an induced subgraph of $G$, we see that 
\begin{equation*}
\text{box}(G)\geq \text{box}(\overline{H_2})=\text{box}(\overline{H_1})+\text{box}(\overline{H_2})
\end{equation*} 
holds.\ In what follows, we may assume that $H_1$ and $H_2$ are nontrivial.\
 
The assumption $d_{\overline{G}}(H_1,H_2)\geq 2$ means that $d_{\overline{G}}(v_1,v_2)\geq 2$ 
for a vertex $v_1$ of $H_1$ and a vertex $v_2$ of $H_2$.\ Hence, an edge of $H_1$ and an edge of 
$H_2$ form $2K_2$, the disjoint union of two edges, as an induced subgraph of $\overline{G}$.\ 
Moreover, we claim the following. \vspace{0.2cm}\\ 
{\bf Claim (1)}: no cointerval spanning subgraphs of $\overline{G}$ contain an edge of $H_1$ 
and an edge of $H_2$, and  \\ 
{\bf Claim (2)}: we need at least $\text{box}(\overline{H_i})$ cointerval spanning subgraphs of 
$\overline{G}$ to cover all edges of $H_i$, where $i=1,2$.\vspace{0.2cm}\\ 
Claim (1) follows from the forbidden subgraph characterization of cointerval graphs.\ 
Actually, cointerval graphs do not contain $2K_2$ as an induced subgraph.\ Claim (2) follows from 
Theorem 2.1.\ A cointerval graph with edges of $H_1$ does not contain edges of $H_2$.\  
Thus, the inequality $\text{box}(G)\geq \text{box}(\overline{H_1})+\text{box}(\overline{H_2})$ 
holds.\
\end{proof}
\end{lem}
We can derive a positive answer to Question 2 by using Lemma 3.2.\ The following lemma is useful 
to make our answer more precise.\ Here, $\lceil x \rceil $ denotes the smallest integer at least $x$.\

\begin{lem}[\cite{CR83}, Lemma 3]
Let $G$ be a graph.\ Let $S_1=\{a_1, a_2, \ldots , a_n\}$ and $S_2=\{b_1,b_2, \ldots , b_n\}$ be 
disjoint subsets of $V(G)$ such that the only edges between $S_1$ and $S_2$ in $\overline{G}$ 
are the edges $a_ib_i$, where $i=1,2,\ldots , n$.\ Then, $\text{box}(G)\geq \lceil \frac{n}{2} \rceil $.\   
\end{lem}
\begin{thm}
For a graph $G$ with $l$ universal vertices, the following inequality holds:
\begin{equation*}
\text{box}(M(G)) \geq \text{box}(G)+\left \lceil \frac{l}{2}\right \rceil .
\end{equation*}
\begin{proof}
Let $G$ be a graph and $x_1, x_2, \ldots , x_l$ universal vertices of $G$.\ Let $H$ be the subgraph of $G$ 
induced by $V(G)\setminus \{x_1, x_2, \ldots , x_l\}$.\ Note that $\text{box}(H)=\text{box}(G)$ holds.\ We 
consider the Mycielski graph $M(G)$ and its complement $\overline{M(G)}$.\ 
Let $X_j=\{(x_1)_j, \, (x_2)_j, \ldots , (x_l)_j\}$, the set of vertices in $V(G)_j$ corresponding to 
universal vertices of $G$.\ 
Let $D_l$ be the subgraph of $M(G)$ induced by the union of $X_1$ and $X_2$.\ 
Note that $X_1$ and $X_2$ are disjoint by their definition.\ It is not difficult to check that 
the only edges between $X_1$ and $X_2$ in $\overline{D_l}$ are the edges $(x_i)_1(x_i)_2$, where 
$i=1,2,\ldots , l$.\ Actually, the vertex $(x_i)_1 \in X_1$ is adjacent to all vertices in 
$V(G)_2\setminus \{(x_i)_2\}$ in $M(G)$ and the vertex $(x_i)_2 \in X_2$ is 
adjacent to all vertices in $V(G)_1\setminus \{(x_i)_1\}$ in $M(G)$ since $x_i$ 
is a universal vertex of $G$.\ We see that $\text{box}(D_l)\geq \lceil \frac{l}{2}\rceil $ 
by Lemma 3.3.\ 

We prove that $d_{\overline{M(G)}}(\overline{H}, \overline{D_l})\geq 2$ holds.\ Let $v$ be a vertex of 
$\overline{H}$ and $x$ a vertex of $\overline{D_l}$.\ The vertex $v$ is in $V(G)_1\setminus X_1$ and the 
vertex $x$ is in $X_1$ or $X_2$.\ We may represent $x$ as $(x_i)_j$, where $j=1,2$.\ Since $x_i$ is a universal 
vertex of $G$, the vertex $(x_i)_j$ is not adjacent to $v$ in $\overline{M(G)}$.\ 
This implies that $d_{\overline{M(G)}}(v, x)\geq 2$ for a vertex $v$ of $\overline{H}$ and a vertex 
$x$ of $\overline{D_l}$, that is, $d_{\overline{M(G)}}(\overline{H}, \overline{D_l})\geq 2$.\ 
Thus, the inequality  
\begin{equation*}
\text{box}(M(G)) \geq \text{box}(H)+\text{box}(D_l) \geq \text{box}(G)+\left \lceil \frac{l}{2}\right \rceil 
\end{equation*}
holds by Lemma 3.2.\
\end{proof}
\end{thm}
\begin{rem}
We note the proof of Theorem 3.4 works on the generalized Mycielski graph $M_r(G)$ (see section 5 for definition), 
that is, $\text{box}(M_r(G)) \geq \text{box}(G)+\left \lceil \frac{l}{2}\right \rceil $ holds 
for a graph with $l$ universal vertices.\ Further observations on $\text{box}(M_r(G))$ appear in section 5.\  

In the proof of Theorem 3.4, we prove that $\text{box}(D_l)\geq \lceil \frac{l}{2}\rceil $ 
by using Lemma 3.3.\ Actually, note that $\text{box}(D_l)=\lceil \frac{l}{2}\rceil $.\  
Any two vertices in $X_1$ are not adjacent in $\overline{M(G)}$ since they are 
adjacent in $M(G)$.\ Hence, $X_1$ is independent in $\overline{D_l}$.\ Also note that $X_2$ is a 
clique in $\overline{M(G)}$ by the definition of Mycielski graphs, that is, in $\overline{D_l}$.\ 
Also see the argument behind the proof of Theorem 5 in \cite{CR83}.\ 

If we restrict our attention to the graph $G$ with only one universal vertex or only two universal 
vertices in the proof of Theorem 3.4, then Lemma 3.3 is superfluous.\ Note that $\text{box}(D_1)=
\text{box}(D_2)=1$ since $D_1$ is the trivial graph with two vertices and $D_2$ is the path with 
four vertices.\  
\end{rem}
Theorem 3.4 implies that for a graph $G$ with universal vertices, $\text{box}(M(G))>\text{box}(G)$ holds.\ 
Also note that Mycielski's graph operation can be used to construct graphs with arbitrary large boxicity 
(and chromatic number) the same as the join of graphs.\

At the end of this section, we note that there is a graph $G$ without universal vertices such that 
the boxicity of the Mycielski graph $M(G)$ is more than that of $G$.\ We give a simple example here.\ 
Also see section 6.\ 
\begin{exmp}
Let $P_n$ be a path with $n$ vertices, where $n\geq 2$.\ We see that $\text{box}(M(P_n))=2$.\ 
We can give a representation of $M(P_n)$ by a family of boxes in Euclidean 2-space.\ See Fig.\ 3 below, 
where we write $V(P_n)_1=\{1,2, \ldots , n\}$ and $V(P_n)_2=\{1',2',\ldots , n'\}$ and for a vertex 
$v\in V(M(P_n))=\{z\}\cup V(P_n)_1\cup V(P_n)_2$, $B_v$ denotes a box in Euclidean 2-space corresponding 
to the vertex $v$.\ Also note that $M(P_n)$ contains an induced cycle $C_5$.\ 
\begin{figure}[!h]
\centering
\includegraphics[scale=1,clip]{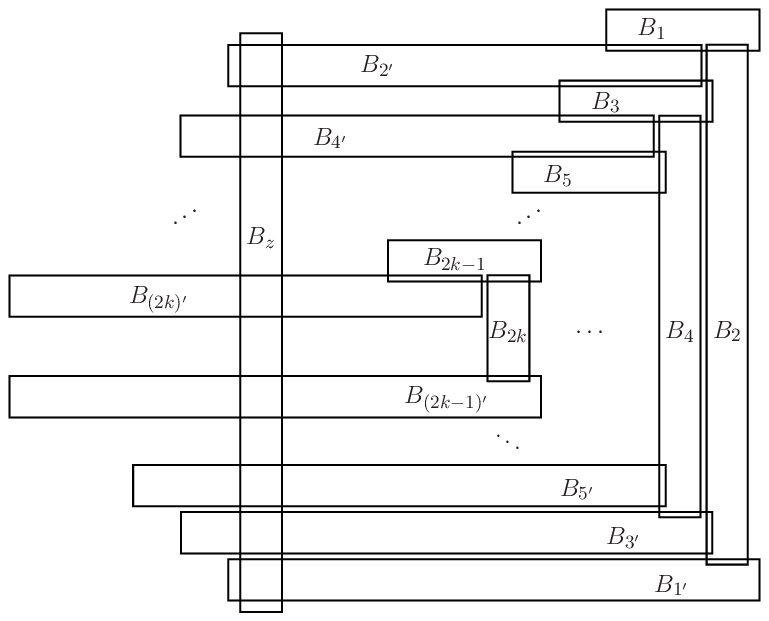}
\renewcommand{\baselinestretch}{1.3}
\caption[]{A representation of $M(P_{2k})$ by a family of boxes in Euclidean 2-space.}
\end{figure}
\end{exmp}
\section{An upper bound for the boxicity of Mycielski graphs}
In this section, we give an upper bound for the boxicity of Mycielski graphs.\
Moreover we calculate the boxicity of Mycielski graphs of some of complete multi-partite graphs.\ 
First we determine the boxicity of Mycielski graphs of complete graphs.\ 

\begin{lem}
For a complete graph $K_n$, the following equalities hold:
\begin{equation*}
\text{box}(M(K_n))=
\begin{cases}
\left\lceil \frac{n}{2}\right\rceil & \text{if $n$ is odd},\\
\left\lceil \frac{n}{2}\right\rceil +1& \text{if $n$ is even}.
\end{cases}
\end{equation*}
\begin{proof}
Let $H_0$ be the subgraph of $M(K_n)$ induced by $V(M(K_n))-\{z\}$.\ We have the inequality 
$\text{box}(M(K_n))\geq \text{box}(H_0)\geq \left\lceil \frac{n}{2}\right\rceil $ by Lemma 3.3.\

Let $V(K_n)=\{v_1,v_2, \ldots , v_n\}$.\ To see $\text{box}(M(K_n))\leq \left\lceil \frac{n}{2}\right\rceil +1$, 
we give cointerval subgraphs of $\overline{M(K_n)}$.\ Let $G_0$ be the subgraph of $\overline{M(K_n)}$ 
induced by $\{z, (v_n)_2\}\cup V(K_n)_1$.\ We define $G_i$ to be the subgraph of $\overline{M(K_n)}$ 
induced by $\{(v_{2i-1})_1, (v_{2i})_1\}\cup V(K_n)_2$, where $i=1,2,\ldots , \lceil \frac{n}{2}\rceil -1$.\ 
Moreover, let $G_{\lceil \frac{n}{2}\rceil }$ be the subgraph of $\overline{M(K_n)}$ induced by 
$\{(v_{n-1})_1, (v_{n})_1\}\cup V(K_n)_2$.\ It is easy to see that the family 
$\{G_0, G_1, \ldots , G_{\lceil \frac{n}{2}\rceil }\}$ is a cointerval edge covering of 
$\overline{M(K_n)}$, and hence $\text{box}(M(K_n))\leq \left\lceil \frac{n}{2}\right\rceil +1$ holds.\ 

If $n$ is odd, the family $\{G_0, G_1, \ldots , G_{\lceil \frac{n}{2}\rceil -1}\}$ is a cointerval edge covering of 
$\overline{M(K_n)}$, because the edge $(v_n)_1(v_n)_2$ is covered with the graph $G_0$.\ Hence we have the equality $\text{box}(M(K_n))=\left\lceil \frac{n}{2}\right\rceil $.\

If $n$ is even, that is, $n=2k$, we show that $\text{box}(M(K_{2k}))>k$.\ Suppose to the contrary that 
$\overline{M(K_{2k})}$ can be covered with $k$ cointerval (spanning) subgraphs $H_1, H_2, \ldots , H_k$ of 
$\overline{M(K_{2k})}$.\ Let $e_j=(v_j)_1(v_j)_2$ for $j=1,2, \ldots , 2k$.\
The graph $H_i$ contains at most two edges in ${\mathcal E}=\{e_1, e_2, \ldots , e_{2k}\}$ since $H_i$ is cointerval.\ 
Actually, the graph $H_i$ must contain two edges in ${\mathcal E}$.\ Otherwise there is a graph 
$H$ in ${\mathcal H}=\{H_1, H_2, \ldots , H_k\}$ which contains only one edge 
in ${\mathcal E}$ or which contains no edges in ${\mathcal E}$.\ 
Hence the family ${\mathcal H} \setminus \{H\}$ of $k-1$ cointerval subgraphs of $\overline{M(K_{2k})}$ must cover 
at least $2k-1$ edges in ${\mathcal E}$, but this is impossible.\ 
On the other hand, there is a cointerval graph $H_*$ in ${\mathcal H}$ which contains an edge $z(v)_1$, where 
the vertex $v$ is in $V(K_{2k})$.\ We may assume that the graph $H_*$ contains two edges $e_s$ and $e_t$ 
in ${\mathcal E}$.\ Hence we see $V(H_*)\supset \{(v_s)_1, (v_s)_2, (v_t)_1, (v_t)_2, z\}$.\ We note that 
\begin{equation*}
(v_s)_1(v_t)_1, (v_s)_1(v_t)_2, (v_t)_1(v_s)_2, z(v_s)_2, z(v_t)_2\not\in E(\overline{M(K_{2k})}) 
\end{equation*} 
by the definition of Mycielski's construction.\ 
\begin{figure}[!h]
\centering
\includegraphics[scale=1,clip]{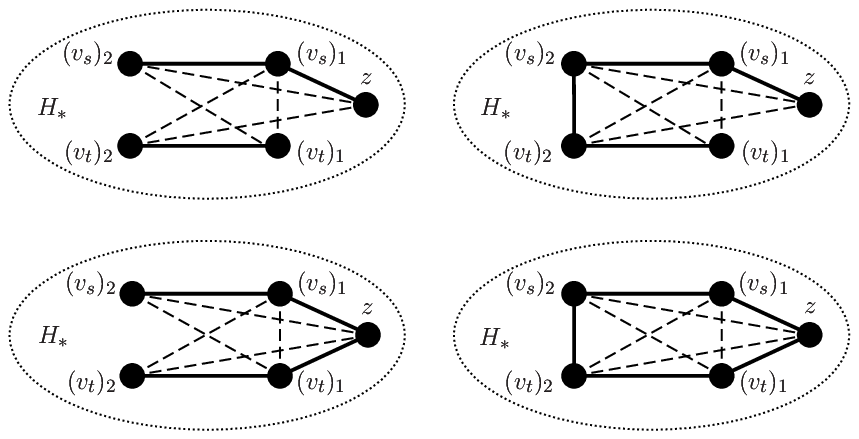}
\renewcommand{\baselinestretch}{1.2}
\caption[]{The subgraph $H_*$ of $\overline{M(K_{2k})}$ containing edges $e_s$ and $e_t$.}
\end{figure}
If $v\not\in \{v_s, v_t\}$, it follows from Lemma 3.3 that 
$\text{box}(\overline{H_*})\geq 2$ since $(v)_1(v_s)_1, (v)_1(v_t)_1\not\in E(\overline{M(K_{2k})})$, 
a contradiction.\ 
Hence we may assume that $v=v_s$.\ We reach the four cases on the graph $H_*$ indicated in Fig.\ 4.\ 
These cases imply that $\text{box}(\overline{H_*})\geq 2$, which contradicts our assumption that $H_*$ is cointerval.\ 
Thus we have $\text{box}(M(K_{2k}))>k$.\ Hence we obtain the equality 
$\text{box}(M(K_n))=\left\lceil \frac{n}{2}\right\rceil +1$ if $n$ is even.
\end{proof}
\end{lem}
\begin{rem}
We proved that the inequality $\text{box}(M(K_n))\leq \left\lceil \frac{n}{2}\right\rceil +1$ holds at 
the second paragraph of the proof of Lemma 4.1.\ We can also derive this inequality by 
using the minimum cardinality of a vertex cover of $M(K_n)$, that is, 
using the inequality $\text{box}(M(K_n))\leq \lfloor \frac{t(M(K_n))}{2}\rfloor +1$.\ 
A subset $U$ of the vertex set of a graph $G$ is a {\it vertex cover} of $G$ if for each $e\in E(G)$, there is a vertex 
$u\in U$ such that $u$ is in $e$.\ Note that $t(M(K_n))=n+1$.\ 
\end{rem} 
The {\it edge clique cover number} of a graph $G$, denoted by $\theta (G)$, is the minimum cardinality of 
a family of cliques that covers all edges of $G$.\ The following theorem gives us an upper bound for the boxicity 
of Mycielski graphs.\ 
\begin{thm}
For a graph $G$ with $l$ universal vertices, the inequality
\begin{equation*}
\text{box}(M(G))\leq \theta (\overline{G})+\left\lceil \frac{l}{2}\right\rceil +1
\end{equation*}
holds.\ If $l$ is zero or odd, we have the inequality 
\begin{equation*}
\text{box}(M(G))\leq \theta (\overline{G})+\left\lceil \frac{l}{2}\right\rceil .
\end{equation*}
\begin{proof}
Let $\{A_1, A_2, \ldots , A_{\theta (\overline{G})}\}$ be a family of cliques in $\overline{G}$ that covers 
all edges of $\overline{G}$.\ Let $v_1, v_2, \ldots , v_l$ be all isolated vertices of $\overline{G}$ and write 
$J=\{v_1, v_2, \ldots , v_l\}$.\ Note that $V(G)=A_1\cup A_2\cup \ldots \cup A_{\theta (\overline{G})} \cup J$.\ 
We define $H_i$ to be the subgraph of $\overline{M(G)}$ induced by $(A_i)_1\cup V(G)_2 \cup \{z\}$ and 
let $E_i=\{xy\,|\, x,y\in V(G)_2\setminus (A_i)_2\}$ and $F_i=\{xy\,|\, x\in (A_i)_1,y\in V(G)_2\setminus (A_i)_2\}$,
where $i=1,2,\ldots , \theta (\overline{G})$.\ We can check that $H_i-(E_i\cup F_i)$ is a cointerval graph 
(see Fig.\ 5).\
\begin{figure}[!h]
\centering
\includegraphics[scale=1,clip]{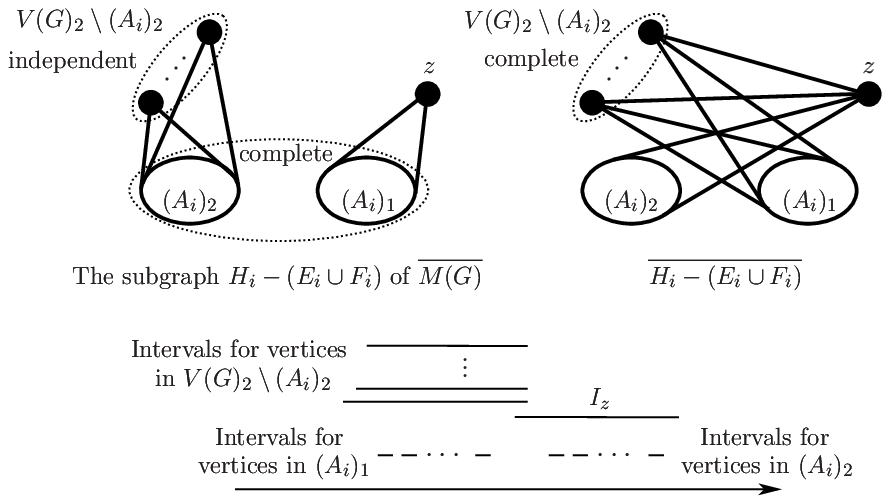}
\renewcommand{\baselinestretch}{1.3}
\caption[]{The subgraph $H_i-(E_i\cup F_i)$ and an interval representation of $\overline{H_i-(E_i\cup F_i)}$.}
\end{figure}
Note that the subgraph of $M(G)$ induced by $J_1\cup J_2\cup \{z\}$ is isomorphic to $M(K_l)$.\ 
Hence the edge set of the subgraph of $\overline{M(G)}$ isomorphic to $\overline{M(K_l)}$ can be covered 
with at most $\left\lceil \frac{l}{2}\right\rceil +1 $ cointerval subgraphs as in the proof of Lemma 4.1.\ 
Let $G_0$ be the subgraph of $\overline{M(G)}$ induced by $\{z, (v_l)_2\}\cup J_1$ and 
$G_i$ the subgraph of $\overline{M(G)}$ induced by $\{(v_{2i-1})_1, (v_{2i})_1\}\cup J_2$ for 
$i=1,2,\ldots , \lceil \frac{l}{2}\rceil -1$.\ Moreover, let $G_{\lceil \frac{l}{2}\rceil }$ be the subgraph of 
$\overline{M(G)}$ induced by $\{(v_{l-1})_1, (v_{l})_1\}\cup J_2$.\
We can check that $\theta (\overline{G})+\left\lceil \frac{l}{2}\right\rceil +1$ cointerval subgraphs 
$H_1-(E_1\cup F_1), \ldots , H_{\theta (\overline{G})}-(E_{\theta (\overline{G})}\cup F_{\theta (\overline{G})}), 
G_0, G_1, \ldots , G_{\lceil \frac{l}{2}\rceil }$ cover all edges of $\overline{M(G)}$.\ 

Let $e$ be an edge of $E(\overline{M(G)})$.\ If $e\cap \{z\}\ne \emptyset $, we see $e\cap V(G)_1\ne \emptyset $.\ 
Hence there is an $i\in \{1,2,\ldots , \theta (\overline{G})\}$ such that $e\in E(H_i-(E_i\cup F_i))$ or $e\in E(G_0)$.\ 
If $e\cap \{z\}=\emptyset $, we have $e\subset V(G)_1\cup V(G)_2$.\ Hence, if $e\subset V(G)_2$, especially, 
$e\cap (A_i)_2\ne \emptyset $, we see $e\in E(H_i-(E_i\cup F_i))$.\ If $e\subset V(G)_2$ and $e\cap (A_i)_2=\emptyset $ 
for any $i$, we see $e\subset J_2$, and hence $e\in E(G_i)$ for $i\ne 0$.\ 
If $e\cap V(G)_1\ne \emptyset $, we reach the following two cases: 
\begin{center}
(i) $e\subset V(G)_1$ or (ii) $e\cap V(G)_2\ne \emptyset $.\  
\end{center}
In the case (i), the edge $e$ is in some $(A_i)_1$ since the family $\{A_1, A_2, \ldots , A_{\theta (\overline{G})}\}$ 
of cliques covers all edges of $\overline{G}$, and hence we have $e\in E(H_i-(E_i\cup F_i))$.\ 

Now we focus on the case (ii).\ Let $u$ be a vertex in $V(\overline{G})$ and $C_u$ the union of cliques in 
$\{A_1, A_2, \ldots , A_{\theta (\overline{G})}\}$ containing the vertex $u$.\ If $u$ is an isolated vertex in 
$\overline{G}$, let $C_u$ be the set $\{u\}$.\ Then we note $u_1\in V(G)_1$ is never adjacent to vertices in 
$V(G)_2\setminus (C_u)_2$ on $\overline{M(G)}$ by the definition of Mycielski graphs.\ Hence the following two 
cases occur: 
\begin{quote}
(ii-1) the edge $e$ connects a vertex of $(A_i)_1$ and a vertex of $(A_i)_2$ for some $i$ or\\ 
(ii-2) the edge $e$ connects a vertex $(v_i)_1$ and a vertex $(v_i)_2$, where $v_i\in J$.\ 
\end{quote}
Under the case (ii-1), we notice $e\in E(H_i-(E_i\cup F_i))$.\ Under the case (ii-2), we see 
$e\in E(G_{\lceil \frac{i}{2}\rceil })$.\ These arguments complete the proof of our first statement.\  

If $l=0$, the graphs 
$H_1-(E_1\cup F_1), \ldots , H_{\theta (\overline{G})}-(E_{\theta (\overline{G})}\cup F_{\theta (\overline{G})})$ 
cover all edges of $\overline{M(G)}$.\ If $l$ is odd, 
$H_1-(E_1\cup F_1), \ldots , H_{\theta (\overline{G})}-(E_{\theta (\overline{G})}\cup F_{\theta (\overline{G})}), 
G_0, G_1, \ldots , G_{\lceil \frac{l}{2}\rceil -1}$ cover all edges of $\overline{M(G)}$, because the edge 
$(v_l)_1(v_l)_2$ is covered with the graph $G_0$.\ Our second statement follows from similar arguments as above.\ 
\end{proof}
\end{thm}
Theorem 3.4 and Theorem 4.3 pretty much narrow the boxicity of Mycielskians of graphs that 
satisfy the equality $\text{box}(G)=\theta (\overline{G})$.\ They also determine the boxicity of some Mycielski graphs.\ 
\begin{cor}
For a graph $G$ with $l$ universal vertices that satisfies the equality $\text{box}(G)=\theta (\overline{G})$, 
the inequalities    
\begin{equation*}
\text{box}(G)+ \left \lceil \frac{l}{2}\right \rceil \leq \text{box}(M(G))
\leq \text{box}(G)+ \left \lceil \frac{l}{2}\right \rceil +1
\end{equation*}
hold.\ Moreover if $l$ be zero or odd, the equality
\begin{equation*}
\text{box}(M(G))=\text{box}(G)+ \left \lceil \frac{l}{2}\right \rceil  
\end{equation*}
holds. \qed 
\end{cor}
We can give examples of graphs that satisfy $\text{box}(G)=\theta (\overline{G})$.\ 
Recall that the boxicity of a complete $k$-partite graph $K_{n_1, n_2, \ldots, n_k}$ is the number of $n_i$ 
which is at least 2.\ If $K_{n_1, n_2, \ldots, n_k}$ has $l$ universal vertices, 
we obtain $\text{box}(K_{n_1, n_2, \ldots, n_k})=k-l=\theta (\overline{K_{n_1, n_2, \ldots, n_k}})$.\ 
Hence we have   
\begin{equation*}
\text{box}(K_{n_1,n_2,\ldots , n_k})+\left \lceil \frac{l}{2}\right \rceil 
=\left \lceil \frac{2k-l}{2}\right \rceil .
\end{equation*}
\begin{cor}
For a complete $k$-partite graph $K_{n_1, n_2, \ldots, n_k}$ with $l$ universal vertices, the inequalities   
\begin{equation*}
\left \lceil \frac{2k-l}{2}\right \rceil \leq \text{box}(M(K_{n_1,n_2,\ldots , n_k}))\leq 
\min\left \{k, \left \lceil \frac{2k-l}{2}\right \rceil +1\right \}
\end{equation*}
hold.\ Especially, if $l$ is zero or odd, the equality 
$\text{box}(M(K_{n_1,n_2,\ldots , n_k}))=\left \lceil \frac{2k-l}{2}\right \rceil $ holds.\ \qed 
\end{cor}
We present other examples of graphs that satisfy $\text{box}(G)=\theta (\overline{G})$.\ 
The graph $H$ whose complement is a chain of cliques is a desired one, where neighboring cliques 
share exactly one vertex and each clique has at least 4 vertices.\ Note that the graph $H$ contains a complete 
multi-partite graph $K_{2,2,\ldots ,2}$ as an induced subgraph and the number of 
its partite sets is equal to that of maximal cliques of the complement $\overline{H}$.\ 

Moreover if we consider a graph operation that extends a graph $G$ to a new graph $\text{Susp}_n(G)$, called the 
{\it $n$-suspension} of $G$, we can get more examples that we desire.\ 
The vertex set of $\text{Susp}_n(G)$ is the union of $V(G)$ and the set of new 
vertices $\{x_1, x_2, \ldots , x_n\}$.\ The edge set of $\text{Susp}_n(G)$ is the union of $E(G)$ and the set 
$\{x_iv\,|\, v\in V(G), i=1,2,\ldots , n \}$.\ Here we assume that $n$ is an integer at least 2.\ We see that 
$\text{box}(\text{Susp}_n(G))=\text{box}(G)+1$ and $\theta (\overline{\text{Susp}_n(G)})=\theta (\overline{G})+1$ 
for a graph $G$ by Theorem 2.1 and Lemma 3.2.\ Hence if the graph $G$ satisfies 
$\text{box}(G)=\theta (\overline{G})$, the equality $\text{box}(\text{Susp}_n(G))=\theta (\overline{\text{Susp}_n(G)})$ 
holds.\ We note that the family of graphs satisfying $\text{box}(G)=\theta (\overline{G})$ is not narrow at all.\

\section{Relation between boxicity of Mycielski graphs and generalized Mycielski graphs}
In this section, we consider relations between Mycielski graphs and their analogous ones 
in the context of boxicity.\  

Let $G$ be a graph and $r$ an integer at least 2.\ Let $V(G)_i$ be a copy of $V(G)$, where 
$i=1,2,\ldots , r$.\ For each vertex $v\in V(G)$, the symbol $v_i$ denotes the vertex in $V(G)_i$ corresponding 
to $v$.\ The {\it generalized Mycielski graph} of $G$, denoted by $M_r(G)$, is the graph 
whose vertex set is $\displaystyle{ \{z\} \cup \bigcup _{i=1}^r V(G)_i }$, the disjoint union of the set of 
an additional new vertex $z$ and copies $V(G)_1, \ldots , V(G)_r$ of $V(G)$, and whose edge set is 
$\displaystyle{\bigcup _{i=1}^{r+1} E_i }$, where
\begin{align*}
E_1&=\{u_1v_1\,|\,uv\in E(G)\},\\
E_i&=\{u_{i-1}v_i, v_{i-1}u_i\, |\, uv\in E(G)\} \hspace{0.5cm} \text{for}\,\,i=2,3,\ldots , r, \text{and}\\
E_{r+1}&=\{zu_r \, |\,u \in V(G)\}.
\end{align*}
Note that the graph $M_2(G)$ is identical to $M(G)$.\ First, we present a relation between $\text{box}(M_r(G))$ 
and $\text{box}(M_2(G))$ for a bipartite graph $G$.\ 
\begin{thm}
For a bipartite graph $G$ and $r\geq 2$, the inequality 
$\text{box}(M_r(G))\leq \text{box}(M_2(G))+2$ holds.\
\begin{proof}
We partition $V(G)$ into two partite sets $V_1$ and $V_2$.\ Fix a family $\{B_x\}$ of boxes in the optimal dimensional 
space which represents $M_2(G)$.\ 
Note that $B_{u_1}\cap B_{u_2}=\emptyset $, $B_{v_1}\cap B_{v_2}=\emptyset $, 
and $B_{u_2}\cap B_{v_2}=\emptyset $ for distinct two vertices $u$ and $v$ of $G$ by the definition of $M_2(G)$.\
Moreover we note that $uv\in E(G)$, $B_{u_1}\cap B_{v_1}\ne \emptyset $,  
$B_{u_1}\cap B_{v_2}\ne \emptyset $, and $B_{u_2}\cap B_{v_1}\ne \emptyset $ are equivalent each other.\ 
First we define the family $\{B_{v_i}'\}$ of boxes in $(\text{box}(M_2(G))+1)$-dimensional space 
to give a box-representation of the graph $M_r(G)-\{z\}$ as follows:\ for each vertex $v\in V(G)$, 
\begin{align*}
B_{v_1}'&=B_{v_1}\times \{0\}, \\ 
B_{v_{2i}}'&=
\begin{cases}
B_{v_2}\times [i-1, i-\frac{1}{2}] & \hspace{0.4cm}\text{if $v\in V_1$},\\
B_{v_2}\times [-(i-\frac{1}{2}), -(i-1)]& \hspace{0.4cm}\text{if $v\in V_2$}, 
\end{cases}
\hspace{1cm}\text{for $i\in \{1,2,\ldots , \lfloor \tfrac{r}{2}\rfloor \}$, and }\\
B_{v_{2i-1}}'&=
\begin{cases}
B_{v_1}\times [-(i-1), -(i-\frac{3}{2})] & \hspace{0.4cm}\text{if $v\in V_1$},\\
B_{v_1}\times [i-\frac{3}{2}, i-1]& \hspace{0.4cm}\text{if $v\in V_2$},
\end{cases}
\hspace{1cm}\text{for $i\in \{2,3,\ldots , \lceil \tfrac{r}{2}\rceil \}$.}
\end{align*}
Take a vertex $v\in V(G)$ and $k\in [r]$, and then consider the adjacency of the vertex $v_k$ 
of $M_r(G)-\{z\}$ from the above family $\{B_{v_i}'\}$ we defined.\ It is easy to see that the box 
$B_{v_k}'$ does not have intersection with boxes corresponding to vertices in 
$\{v_1, v_2, \ldots , v_r\}\setminus \{v_k\}$.\ 
We also see that $B_{v_k}'$ does not have intersection with boxes that correspond to vertices in  
$V(G)_k\setminus \{v_k\}$ for $k\in \{2,3,\ldots , r\}$.\ Clearly, the family $\{B_{v_1}'\}_{v\in V(G)}$ 
represents the subgraph of $M_r(G)-\{z\}$ induced by $V(G)_1$, so we may assume $k\geq 2$.\ 

If the vertex $v$ is adjacent to a vertex $u$ in $G$, we can check that the box $B_{v_k}'$ has nonempty 
intersection only with boxes $B_{u_{k-1}}'$ and $B_{u_{k+1}}'$ for $2\leq k\leq r-1$, and only 
with the box $B_{u_{r-1}}'$ for $k=r$ in the family 
$\{B_{u_1}', B_{u_2}', \ldots , B_{u_r}'\}$.\ If the vertex $v$ is not adjacent to a vertex $u$ in $G$, 
no boxes in the family $\{B_{u_1}', B_{u_2}', \ldots , B_{u_r}'\}$ have nonempty intersection 
with $B_{v_k}'$ since $B_{u_1}\cap B_{v_1}= \emptyset $, $B_{u_1}\cap B_{v_2}= \emptyset $, 
and $B_{u_2}\cap B_{v_1}= \emptyset $ hold.\ Hence the family $\{B_{v_i}'\}$ represents the graph $M_r(G)-\{z\}$.\
 
Finally, we define the  family $\{B_x''\}$ of boxes in $(\text{box}(M_2(G))+2)$-dimensional space 
as follows: that represents of $M_r(G)$ as follows:\
\begin{align*}
B_{v_i}''&=B_{v_i}'\times \{0\} \hspace{0.4cm}\text{for $i\ne r$}\\
B_{v_r}''&=B_{v_r}'\times [0, 1]\\
B_{z}''&=B\times \{1\},
\end{align*}
where $B$ is a box in $(\text{box}(M_2(G))+1)$-dimensional space that contains all boxes in $\{B_{v_r}'\}_{v\in V(G)}$.\  
We can check easily that the family $\{B_x''\}$ represents of $M_r(G)$, which completes the proof of our theorem.\ 
\end{proof} 
\end{thm}
The author think that the inequality $\text{box}(M_r(G))\leq \text{box}(M_2(G))+c$ holds for a graph $G$ and some 
small constant $c$ in general.\ The next theorem shows a relation between 
$\text{box}(M_r(G))$ and $\text{box}(M_3(G))$ for a graph $G$.\ 
These results will become our motivation a bit to calculate the boxicity of $M_2(G)$ or $M_3(G)$.\
\begin{thm}
For a graph $G$ and $r\geq 3$, the inequality $\text{box}(M_r(G))\leq \text{box}(M_3(G))+1$ holds.\
\begin{proof}
Let $\{B_x\}$ be a family of boxes in the optimal dimensional space which represents $M_3(G)$.\ 
We note that for distinct two vertices $u$ and $v$ of $G$, 
\begin{align*}
B_{u_i}\cap B_{u_j}&=\emptyset ,\, B_{v_i}\cap B_{v_j}=\emptyset ,\, \text{where $i,j\in \{1,2,3\}$ and $i\ne j$},\\ 
B_{u_i}\cap B_{v_i}&=\emptyset ,\, \text{where $i\in \{2,3\}$, and} \\
B_{u_1}\cap B_{v_3}&=\emptyset ,\, B_{v_1}\cap B_{u_3}=\emptyset 
\end{align*}
hold by the definition of $M_3(G)$.\ In addition we note that $uv\in E(G)$, 
$B_{u_1}\cap B_{v_1}\ne \emptyset $, $B_{u_1}\cap B_{v_2}\ne \emptyset $, $B_{u_2}\cap B_{v_1}\ne \emptyset $, 
$B_{u_2}\cap B_{v_3}\ne \emptyset $, and $B_{u_3}\cap B_{v_2}\ne \emptyset $ are equivalent each other.\
By using similar techniques in the previous theorem, we can present the family $\{B_x'\}$ of boxes in 
$(\text{box}(M_3(G))+1)$-dimensional space that represents the graph $M_r(G)$ as follows: for each vertex $v\in V(G)$, 
define 
\begin{align*}
B_{v_i}'&=B_{v_i}\times \{0\} \hspace{2.48cm}\text{for $i\in \{1,2\}$}, \\ 
B_{v_{2i-1}}'&=B_{v_3}\times [i-2, i-\tfrac{3}{2}] \hspace{1cm}\text{for $i\in \{2,3,\ldots , \lceil \tfrac{r}{2}\rceil \}$}, \\
B_{v_{2i}}'&=B_{v_2}\times [i-\tfrac{3}{2}, i-1] \hspace{1cm}\text{for $i\in \{2,3,\ldots , \lfloor \tfrac{r}{2}\rfloor \}$}, 
\end{align*}
and for the additional vertex $z'$ of $M_r(G)$, 
\begin{equation*}
B_{z'}'=
\begin{cases}
B\times \{\tfrac{r}{2}-1 \} \hspace{1.53cm}\text{if $r$ is even}, \\ 
B_z\times \{\lfloor \tfrac{r}{2}\rfloor -\frac{1}{2} \} \hspace{1cm}\text{if $r$ is odd}, 
\end{cases}
\end{equation*}
where $B$ is a box in $(\text{box}(M_3(G)))$-dimensional space that contains all boxes in $\{B_{v_2}\}_{v\in V(G)}$ and 
$z$ is the additional vertex of $M_3(G)$.\ 
Note that any pair of distinct two boxes in $\{B_{v_1}', B_{v_2}',\ldots , B_{v_r}'\}$ does not have 
intersection for a vertex $v$ of $G$, and also note that any pair of distinct two boxes in $\{B_{v_k}'\}_{v\in V(G)}$ 
does not have intersection for $k\in \{2,3,\ldots , r\}$.\

Fix a vertex $v\in V(G)$ and $k\in [r]$.\ We consider the adjacency of the vertex $v_k$ of $M_r(G)$.\ 
Clearly, the family $\{B_{v_1}'\}_{v\in V(G)}$ represents the subgraph of $M_r(G)$ induced by $V(G)_1$, 
and hence we may assume $k\geq 2$.\
If the vertex $v$ is adjacent to a vertex $u$ in $G$, we can verify that the box $B_{v_k}'$ has nonempty 
intersection only with boxes $B_{u_{k-1}}'$ and $B_{u_{k+1}}'$ for $2\leq k\leq r-1$, and only 
with $B_{u_{r-1}}'$ for $k=r$ in $\{B_{u_1}', B_{u_2}', \ldots , B_{u_r}'\}$.\
If the vertex $v$ is not adjacent to a vertex $u$ in $G$, no boxes in the family 
$\{B_{u_1}', B_{u_2}', \ldots , B_{u_r}'\}$ have nonempty intersection 
with $B_{v_k}'$.\ In addition, the box $B_{v_k}'$ have nonempty intersection with $B_{z'}'$ if and only if $k=r$
for each $v\in V(G)$.\ Hence our arguments guarantee that the family $\{B_x'\}$ represents the graph $M_r(G)$.\
\end{proof}
\end{thm}

\section{Concluding Remarks: graphs with $\text{box}(M(G))>\text{box}(G)$} 
We proved that the boxicity of the Mycielski graph of a graph $G$ with universal vertices is 
more than that of $G$.\ As examples of complete multi-partite graphs without universal vertices, 
one may expect that the equality $\text{box}(M(G))=\text{box}(G)$ holds for a graph $G$ without 
universal vertices.\ However, we note that there is a graph $G$ without universal vertices such that 
$\text{box}(M(G))>\text{box}(G)$ holds.\ For examples, nontrivial interval graphs without universal 
vertices are the desired ones.\ The Mycielski graph of such an interval graph is not interval because it 
contains a cycle with 5 vertices as an induced subgraph.\ Another example of a graph without universal 
vertices that satisfy $\text{box}(M(G))>\text{box}(G)$ is a cycle with at least 5 vertices.\ 
The author verified in manuscript that the boxicity of Mycielski graph of a cycle with at 
least 5 vertices is equal to 3.\

\section*{Acknowledgments}
This work was supported by Grant-in-Aid for Young Scientists (B), No.25800091.

\end{document}